\documentclass{article}%
\usepackage{graphicx}
\usepackage{amsmath}
\usepackage{amsfonts}
\usepackage{amssymb}%
\newtheorem{theorem}{Theorem}
{}

\newtheorem{lemma}{Lemma}
{}

\newenvironment{proof}[1][Proof]{\textbf{#1.} }{\ \rule{0.5em}{0.5em}}

\begin{document}

\title{On\ the Hill's Operator with Matrix Potential}
\author{O.A.Veliev\\{\small Depart. of Math, Fac.of Arts and Sci.,}\\{\small Dogus University, Ac\i badem, 34722,}\\{\small Kadik\"{o}y, \ Istanbul, Turkey.}\\\ {\small e-mail: oveliev@dogus.edu.tr}}
\date{}
\maketitle

\begin{abstract}
In this article we obtain the asymptotic formulas for the eigenvalues and
eigenfunctions of the self-adjoint operator generated by a system of
Sturm-Liouville equations with summable coefficients and the quasiperiodic
boundary conditions. Then using these asymptotic formulas, we find the
conditions on the potential for which the number of gaps in the spectrum of
the Hill's operator with matrix potential is finite.

\end{abstract}

\bigskip Let $L(Q(x))$ be the differential operator generated in the space
$L_{2}^{m}(-\infty,\infty)$ of the vector functions by the differential
expression
\begin{equation}
-y^{^{\prime\prime}}(x)+Q\left(  x\right)  y(x),
\end{equation}
where $Q\left(  x\right)  =\left(  b_{i,j}\left(  x\right)  \right)  $ is a
$m\times m$ Hermitian matrix with the complex-valued summable entries
$b_{i,j}\left(  x\right)  $ and $Q\left(  x+1\right)  =Q\left(  x\right)  $.
It is well-known that the spectrum of the operator $L$ is the union of the
spectra of the operators $L_{t}$ for $t\in\lbrack0,2\pi)$ generated in
$L_{2}^{m}(0,1)$ by the differential expression (1) and the quasiperiodic
conditions
\begin{equation}
y^{^{\prime}}\left(  1\right)  =e^{it}y^{^{\prime}}\left(  0\right)  ,\text{
}y\left(  1\right)  =e^{it}y\left(  0\right)  .
\end{equation}
Note that $L_{2}^{m}(a,b)$ is the set of the vector functions

$f\left(  x\right)  =\left(  f_{1}\left(  x\right)  ,f_{2}\left(  x\right)
,...,f_{m}\left(  x\right)  \right)  $ with $f_{k}\left(  x\right)  \in
L_{2}(a,b)$ for $k=1,2,...,m.$ The norm $\left\Vert .\right\Vert $ and inner
product $(.,.)$ in $L_{2}^{m}(a,b)$ are defined by%
\[
\left\Vert f\right\Vert =\left(  \int\limits_{a}^{b}\left\vert f\left(
x\right)  \right\vert ^{2}dx\right)  ^{\frac{1}{2}},\text{ }(f\left(
x\right)  ,g\left(  x\right)  )=\int\limits_{a}^{b}\left\langle f\left(
x\right)  ,g\left(  x\right)  \right\rangle dx,
\]
where $\left\vert .\right\vert $ and $\left\langle .,.\right\rangle $ are the
norm and inner product in $\mathbb{C}^{m}.$

Let us introduce some preliminary results and describe briefly the scheme of
the paper. Clearly,

$\varphi_{k,1,t}=\left(
\begin{array}
[c]{c}%
e^{i\left(  2\pi k+t\right)  x}\\
0\\
\vdots\\
0
\end{array}
\right)  ,\varphi_{k,2,t}=\left(
\begin{array}
[c]{c}%
0\\
e^{i\left(  2\pi k+t\right)  x}\\
\vdots\\
0
\end{array}
\right)  ,...,\varphi_{k,m,t}=\left(
\begin{array}
[c]{c}%
0\\
\vdots\\
0\\
e^{i\left(  2\pi k+t\right)  x}%
\end{array}
\right)  $ are the eigenfunctions of the operator $L_{t}(0)$ corresponding to
the eigenvalue $\left(  2\pi k+t\right)  ^{2}$. If $t\neq0,\pi$ then the
multiplicity of the eigenvalue $\left(  2\pi k+t\right)  ^{2}$ is $m$ and the
corresponding eigenspace is $E_{k}(t)=span\left\{  \varphi_{k,1,t}%
,\varphi_{k,2,t},...,\varphi_{k,m,t}\right\}  $. In the cases $t=0$ and
$t=\pi$ the multiplicity of the eigenvalues $\left(  2\pi k\right)  ^{2}$ and
$\left(  2\pi k+\pi\right)  ^{2}$ is $2m$ and the corresponding eigenspaces are

$E_{k}(0)=span\left\{  \varphi_{n,j,0}:n=k,-k;\text{ }j=1,2,...m\right\}  $ and

$E_{k}(\pi)=span\left\{  \varphi_{n,j,\pi}:n=k,-(k+1);\text{ }%
j=1,2,...m\right\}  $ respectively

It easily follows from the classical investigations [1] that the eigenvalues
of the operator $L_{t}\left(  Q\right)  $ consist of $m$ sequences
\begin{equation}
\{\lambda_{k,1}(t):k\in\mathbb{Z}\mathbb{\}},\text{ }\{\lambda_{k,2}%
(t):k\in\mathbb{Z}\mathbb{\}},...,\text{ }\{\lambda_{k,m}(t):k\in
\mathbb{Z}\mathbb{\}}%
\end{equation}
satisfying the following, uniform with respect to $t\in\lbrack0,2\pi),$
asymptotic formulas
\begin{equation}
\lambda_{k,j}(t)=\left(  2\pi k+t\right)  ^{2}+O\left(  k^{1-\frac{1}{2m}%
}\right)
\end{equation}
for $j=1,2,...,m.$ We say that the formula \ $f(k,t)=O(h(k))$ is uniform with
respect to $t\in\lbrack0,2\pi)$ if there exists a positive constant $c$,
independent on $t,$ such that $\mid f(k,t))\mid<c\mid h(k)\mid$ for all
$t\in\lbrack0,2\pi)$ and $k\in\mathbb{Z}$.

In forthcoming relations by $N$ we denote a big positive integer, that is,
$N\gg1$, and by $c_{k}$ for $k=1,2,...$, the positive constants, independent
on $N$ and $t,$ whose exact values are inessential. The formula (4) shows that
the eigenvalue $\lambda_{k,j}(t)$ of $L_{t}(Q)$ is close to the eigenvalue
$\left(  2k\pi+t\right)  ^{2}$ of $L_{t}(Q),$ namely :
\[
|\lambda_{k,j}(t)-\left(  2k\pi+t\right)  ^{2}|<c_{1}|k|^{1-\frac{1}{2m}}.
\]
To analyze the distance of the eigenvalue $\left(  2k\pi+t\right)  ^{2}$ of
$L_{t}(0)$ from the other eigenvalue $\left(  2p\pi+t\right)  ^{2}$ of
$L_{t}(Q)$ , which is important in perturbation theory, we have to consider
three cases.

Case 1 : $t\in T(k),$ where

$T(k)=[-\frac{\pi}{2},\frac{3\pi}{2})\backslash((-(\ln|k|)^{-1},(\ln
|k|)^{-1})\cup(\pi-(\ln|k|)^{-1},\pi+(\ln|k|)^{-1}))$ ,

Case 2: $t\in((-(\ln|k|)^{-1},(\ln|k|)^{-1}),$

Case 3: $t\in(\pi-(\ln|k|)^{-1},\pi+(\ln|k|)^{-1})).$

Using (4) one can easily verify that the inequalities
\begin{align}
|\left(  2k\pi+t\right)  ^{2}-\left(  2\pi p+t\right)  ^{2}| &  >c_{2}%
(\ln|k|)^{-1}(||k|-|p||+1)(|k|+|p|),\nonumber\\
|\lambda_{k,j}(t)-\left(  2\pi p+t\right)  ^{2}| &  >c_{2}(\ln|k|)^{-1}%
(||k|-|p||+1)(|k|+|p|)
\end{align}
hold in Case 1, Case 2, and Case 3 for $p\neq k,$ for $p\neq k,-k$ and for

$p\neq k,-(k+1)$ respectively. To avoid the listing of these cases repeatedly
we introduce the following notation%
\begin{equation}
A(k,t)=\left\{
\begin{tabular}
[c]{l}%
$\{k\}$ for $t\in T(k),$\\
$\{k,-k\}$ for $t\in(-\ln|k|)^{-1},(\ln|k|)^{-1}),$\\
$\{k,-k-1\}$ for $t\in(\pi-(\ln|k|)^{-1},\pi+(\ln|k|)^{-1})$%
\end{tabular}
\ \ \ \ \ \ \right\}  .
\end{equation}
Thus the inequalities (5) hold for $p\notin A(k,t).$ This implies the
following purposive relations
\begin{equation}
\sum_{p:p>d}\frac{1}{\left\vert \lambda_{k,j}(t)-(2\pi p+t)^{2}\right\vert
}<\frac{c_{3}}{d},\text{ }\forall d>2\mid k\mid,
\end{equation}%
\begin{equation}
\sum_{p:p\notin A(k,t)}\frac{1}{\left\vert \lambda_{k,j}(t)-(2\pi
p+t)^{2}\right\vert }=O(\frac{\ln|k|}{k}),
\end{equation}%
\begin{equation}
\sum_{p:p\notin A(k,t)}\frac{1}{\left\vert \lambda_{k,j}(t)-(2\pi
p+t)^{2}\right\vert ^{2}}=O(\frac{(\ln|k|)^{2}}{k^{2}})
\end{equation}%
\begin{equation}
\sum_{p:p\notin A(k,t)}\frac{1}{\left\vert \lambda_{p,j}(t)-(2\pi
k+t)^{2}\right\vert ^{2}}=O(\frac{(\ln|k|)^{2}}{k^{2}})
\end{equation}
The formulas (7)-(10) are uniform with respect to $t\in\lbrack-\frac{\pi}%
{2},\frac{3\pi}{2}).$

In this paper we suggest a method by which the asymptotic formulas of high
accuracy for the eigenvalues $\lambda_{k,j}(t)$ and for the corresponding
normalized eigenfunctions $\Psi_{k,j,t}(x)$ of $L_{t}(Q)$ when the entries
$b_{i,j}(x)$ of $Q(x)$ belong to $L_{1}[0,1]$, that is, when there is not any
condition about smoothness of the\ potential $Q(x)$ are obtained. Note that to
obtain the asymptotic formulas of high accuracy by using the classical
asymptotic expansions for the solutions of the matrix equation $-Y^{^{\prime
\prime}}+Q\left(  x\right)  Y=\lambda Y$ it is required \ that $\ Q(x)$ be
differentiable (see [1-4]). To obtain the asymptotic formulas we consider the
operator $L_{t}(Q)$ as perturbation of $L_{t}(C),$ where $C=\int
\limits_{0}^{1}Q\left(  x\right)  dx,$\ by $Q(x)-C,$ that is, we take the
operator$L_{t}(C)$ for an unperturbed operator and the operator of
multiplication by $Q(x)-C$ for a perturbation. Therefore first we analyze the
eigenvalues and eigenfunction of $L_{t}(C)$.

The adjoint operator to $L_{t}\left(  Q\right)  $ is $L_{t}\left(  Q^{\ast
}\right)  ,$ where $Q^{\ast}\left(  x\right)  $ is the adjoint matrix to
$Q\left(  x\right)  $. Since $Q(x)$ is Hermitian matrix, that is, $Q^{\ast
}\left(  x\right)  =Q\left(  x\right)  $ and the boundary conditions (2) are
self-adjoint the operators $L_{t}\left(  Q\right)  ,$ $L_{t}(C)$ and
$L_{t}(0)$ are self-adjoint. The eigenvalues of $C,$ counted with
multiplicity, and the corresponding orthonotmal eigenvectors are denoted by

$\mu_{1}\leq\mu_{2}\leq...\leq\mu_{m}$ and $v_{1},v_{2},...,v_{m}$
respectively. Thus%
\[
Cv_{j}=\mu_{j}v_{j},\text{ }\left\langle v_{i},v_{j}\right\rangle
=\delta_{i,j}%
\]
In these notations the eigenvalues and eigenfunctions of $L_{t}(C)$ are

$\ \ \ \ \ \ \mu_{k,j}(t)=\left(  2\pi k+t\right)  ^{2}+\mu_{j},$
$\Phi_{k,j,t}(x)=v_{j}e^{i\left(  2\pi k+t\right)  x},$ that is,
\begin{equation}
(L(C)-\mu_{k,j}(t))\Phi_{k,j,t}(x)=0.
\end{equation}

To prove the asymptotic formulas for the eigenvalues $\lambda_{k,j}(t)$ and
for the corresponding normalized eigenfunctions $\Psi_{k,j,t}(x)$ of
$L_{t}(Q)$ we use the formula
\begin{equation}
(\lambda_{k,j}(t)-\mu_{n,i}(t))(\Psi_{k,j,t},\Phi_{n,i,t})=((Q(x)-C)\Psi
_{k,j,t},\Phi_{n,i,t})
\end{equation}
for $n\in A(k,t),$ connecting the eigenvalues and eigenfunctions of the
operators $L_{t}(Q)$ and $L_{t}(C),$ which can be obtained from
\begin{equation}
L(Q\left(  x\right)  )\Psi_{k,j,t}(x)=\lambda_{k,j}(t)\Psi_{k,j,t}(x)
\end{equation}
by multiplying both sides by $\Phi_{n,i,t}(x)$ and using (11). Then we
estimate the right-hand side of (12) ( see Lemma 2) by using Lemma 1. At last
using the connecting formula (12) and lemmas 2, 3 we find the asymptotic
formulas for the eigenvalues and eigenfunctions of $L_{t}(Q)$ ( see Theorem
1). Then using these asymptotic formulas, we find the conditions on the
eigenvalues of the matrix $C$ for which the number of the gaps in the spectrum
of the Hill's operator $L(Q)$ is finite. This result for differentiable $Q(x)$
is obtained in [4,5]. To estimate the right-hand side of (12)\ we use (7),
(8), the following lemma, and the formula
\begin{equation}
\left(  \lambda_{k,j}(t)-\left(  2\pi n+t\right)  ^{2}\right)  \left(
\Psi_{k,j,t},\varphi_{n,s,t}\right)  =\left(  \Psi_{k,j,t},Q\left(  x\right)
\varphi_{n,s,t}\right)  ,
\end{equation}
which can be obtained from (13) by multiplying both sides by $\varphi
_{n,s,t}(x)$ and using $L_{t}\left(  0\right)  \varphi_{n,s,t}(x)=\left(  2\pi
n+t\right)  ^{2}\varphi_{n,s,t}(x)$.

\begin{lemma}
For the right-hand side of (14) the followings hold:
\begin{equation}
\left(  \Psi_{k,j,t}(x),Q\left(  x\right)  \varphi_{n,s,t}(x)\right)
=\sum\limits_{\substack{q=1,2,...m\\p=-\infty,...,\infty}}b_{s,q,n-p}%
(\Psi_{k,j,t},\varphi_{p,q,t}),
\end{equation}%
\begin{equation}
\text{ }\left\vert \left(  \Psi_{k,j,t}(x),Q\left(  x\right)  \varphi
_{n,s,t}(x)\right)  \right\vert <c_{4}\text{ }%
\end{equation}
for $n\in Z$ $;$ $\mid k\mid\geq N$ $;$ $s,j=1,2,...,m$, where $\ $%
\begin{equation}
b_{s,q,n-p}=\int\limits_{0}^{1}b_{s,q}\left(  x\right)  e^{2\pi i\left(
p-n\right)  x}dx.
\end{equation}

\end{lemma}

\begin{proof}
Since $Q(x)\Psi_{k,j,t}(x)\in L_{1}^{m}[0,1]$ we have
\[
\lim_{n\rightarrow\infty}\left(  Q(x)\Psi_{k,j,t}(x),\varphi_{n,s,t}%
(x)\right)  =0,\text{ }\forall s=1,2,...,m.
\]
Therefore there exists a positive constant $C(k,j)$ and indices $k_{0},$
$j_{0}$ satisfying
\begin{equation}
\max_{\substack{n\in\mathbb{Z},\\s=1,2,...,m}}\left\vert (\Psi_{k,j,t}%
,Q\left(  x\right)  \varphi_{n,s,t})\right\vert =\left\vert (\Psi
_{k,j,t},Q\left(  x\right)  \varphi_{k_{0},j_{0},t})\right\vert =C(k,j)
\end{equation}
Using this, (14), (5), we obtain
\begin{equation}
\mid\left(  \Psi_{k,j,t}(x),\varphi_{p,s,t}(x)\right)  \mid\leq\frac
{C(k,j)}{\mid\lambda_{k,j}(t)-\left(  2\pi p+t\right)  ^{2}\mid}%
\end{equation}
for $p\notin A(k,t)$ and $s,j=1,2,...,m.$ This and (7) imply that
\[
\sum_{p:p>d}\mid\left(  \Psi_{k,j,t}(x),\varphi_{p,s,t}(x)\right)  \mid
<\frac{c_{3}C(k,l)}{d},
\]
where $d>2|k|$ and $t\in\lbrack-\frac{\pi}{2},\frac{3\pi}{2}).$ Therefore the
decomposition of $\Psi_{k,j,t}(x)$ by the orthonormal basis $\{\varphi
_{p,q,t}(x)$:$p\in\mathbb{Z},$ $q=1,2,...,m\}$ is of the form
\begin{equation}
\Psi_{k,j,t}(x)=\sum_{\substack{p:|p|\leq d,\\q=1,2,...,m}}\left(
\Psi_{k,j,t}(x),\varphi_{p,q,t}(x)\right)  \varphi_{p,q,t}(x)+g_{d}(x),
\end{equation}
where $\sup_{x\in\lbrack0,1]}|g_{d}(x)|<\frac{c_{3}C(k,l)}{d}$. Putting this
in $\left(  \Psi_{k,j,t}(x),Q\left(  x\right)  \varphi_{n,s,t}(x)\right)  $
and tending $d$ to $\infty$, we obtain (15). \ 

Now we prove (16). Using (20) in $(\Psi_{k,j,t}(x),Q\left(  x\right)
\varphi_{k_{0},j_{0},t}(x))$, isolating the terms with multiplicands
$(\Psi_{k,j,t}(x),\varphi_{p,i,t}(x))$ for $p=k,-k,-(k+1);$

$q=1,2,...,m,$ and tending $d$ to $\infty,$we obtain
\begin{align}
\left(  \Psi_{k,j,t}(x),Q\left(  x\right)  \varphi_{k_{0},j_{0,}t}(x)\right)
&  =\sum\limits_{\substack{p=k,-k,-(k+1);\\q=1,2,...,m}}b_{j_{0},q,k_{0}%
-p}(\Psi_{k,j,t}\left(  x\right)  ,\varphi_{p,q,t}\left(  x\right)
)+\nonumber\\
&  \sum\limits_{\substack{p\neq k,-k,-(k+1)\\q=1,2,...,m}}b_{j_{0},q,k_{0}%
-p}\left(  \Psi_{k,j,t}\left(  x\right)  ,\varphi_{p,q,t}\left(  x\right)
\right)  .
\end{align}
Since
\begin{equation}
\mid b_{j,i,s}\mid\leq\max_{p,q=1,2,...,m}\int\limits_{0}^{1}\mid
b_{p,q}\left(  x\right)  \mid dx<c_{5},\forall j,i,s
\end{equation}
( see (17)) it follows from (19), (8), (6) \ that
\[
\sum\limits_{\substack{p\neq k,-k,-(k+1)\\q=1,2,...,m}}b_{j_{0},q,k_{0}%
-p}\left(  \Psi_{k,j,t},\varphi_{p,q,t}\right)  =O(\frac{\ln|k|}{k})C(k,l)).
\]
Therefore taking into account that the absolute value of the first summation
in the right-hand side of (21) is not greater than $3mc_{5}$ ( see (22)), we
conclude that $\mid C(k,l)\mid<4mc_{5}$ ( see (18)) which means that (16) holds
\end{proof}

\begin{lemma}
Let $\Psi_{k,j,t}(x)$ be any normalized eigenfunctions corresponding to the
eigenvalue $\lambda_{k,j}(t)$ of $L_{t}(Q).$ Then%
\begin{equation}
\left(  \Psi_{k,j,t}(x),(Q\left(  x\right)  -C)\Phi_{n,i,t}(x)\right)
=O(\frac{\ln|k|}{k})+O(b_{k})
\end{equation}
for $n\in A(k,t)$ and $i=1,2,...,m,$ where
\begin{equation}
b_{k}=\max\{\mid b_{i,j,n}\mid:i,j=1,2,...m;\text{ }n=2k,-2k,2k+1,-2k-1\}.
\end{equation}
In particular, for $t\in T(k)$ the formula
\begin{equation}
\left(  \Psi_{k,j,t}(x),(Q\left(  x\right)  -C)\Phi_{k,i,t}(x)\right)
=O(\frac{\ln|k|}{k})
\end{equation}
holds. The formula (23) is uniform with respect to $t\in\lbrack-\frac{\pi}%
{2},\frac{3\pi}{2}).$
\end{lemma}

\begin{proof}
Since $\Phi_{n,s,t}(x)\equiv v_{s}e^{i\left(  2\pi n+t\right)  x},$ to prove
(23) it is enough to show that%
\begin{equation}
\left(  \Psi_{k,j,t}(x),(Q\left(  x\right)  -C)\varphi_{n,s,t}(x)\right)
=O((\frac{\ln|k|}{k})+O(b_{k})
\end{equation}
for $s=1,2,...,m$ and $n\in A(k,t).$ Using the obvious relation
\begin{equation}
\left(  \Psi_{k,j,t}(x),C\varphi_{n,s,t}(x)\right)  =\sum\limits_{q=1,2,...,m}%
b_{s,q,0}(\Psi_{k,j,t}(x),\varphi_{n,q,t}(x))\nonumber
\end{equation}
and (15), we see that
\begin{equation}
\left(  \Psi_{k,j,t},(Q\left(  x\right)  -C)\varphi_{n,s,t}(x)\right)
=\sum\limits_{\substack{p:p\in A(k,t)\backslash n\\q=1,2,...,m}}b_{s,q,n-p}%
\left(  \Psi_{k,j,t}(x),\varphi_{p,q,t}(x)\right)  +
\end{equation}%
\[
\sum\limits_{\substack{p:p\notin A(k,t)\\q=1,2,...,m}}b_{s,q,k-p}\left(
\Psi_{k,j,t}(x),\varphi_{p,q,t}(x)\right)  .
\]
On the other hand, it follows from (19), (18), (16) that
\begin{equation}
\mid\left(  \Psi_{k,j,t}(x),\varphi_{p,q,t}(x)\right)  \mid\leq\frac{c_{4}%
}{\mid\lambda_{k,j}(t)-\left(  2\pi p+t\right)  ^{2}\mid}%
\end{equation}
for $p\notin A(k,t)$ and $q=1,2,...,m.$ Using this, (22), and (8) we see that
the second summation of the right-hand side of (27) is $O((\frac{\ln|k|}{k}).$
Besides it follows from (6), (24) that the first summation of the right-hand
side of (27) is $O(b_{k}),$ since for $n\in A(k,t)$ and $p\in A(k,t)\backslash
n$ we have $n-p\in\{2k,-2k,2k+1,-2k-1\}.$ If $t\in T(k)$ then the first
summation of the right-hand side of (27) is absent, since $A(k,t)=\{k\}$ and
$A(k,t)\backslash n=\emptyset$ for $t\in T(k)$ and $n\in A(k,t)$ (see (6)).
Thus (25) is also proved.
\end{proof}

\begin{lemma}
For each eigenfunction $\Psi_{k,j,t}(x)$ of $L_{t}(Q)$ there exists an
eigenfunction $\Phi_{n,i,t}(x)$ of $L_{t}(C)$ satisfying
\begin{equation}
\left\vert \left(  \Psi_{k,j,t}(x),\Phi_{n,i,t}(x)\right)  \right\vert >c_{6}%
\end{equation}
and for each eigenfunction $\Phi_{k,j,t}(x)$ of $L_{t}(C)$ there exists an
eigenfunction $\Psi_{n,i,t}(x)$ of $L_{t}(Q)$ satisfying
\begin{equation}
\left\vert \left(  \Phi_{k,j,t}(x),\Psi_{n,i,t}(x)\right)  \right\vert >c_{7},
\end{equation}
where $n\in A(k,t);$ $i,j=1,2,...,m$ and $|k|\geq N$ .
\end{lemma}

\begin{proof}
It follows from (28) and (9) that
\begin{equation}
\sum\limits_{q=1,2,...,m}(\sum_{n:n\notin A(k,t)}\left\vert \left(
\Psi_{k,j,t},\varphi_{n,q,t}\right)  \right\vert ^{2})=O(\frac{(\ln|k|)^{2}%
}{k^{2}}).
\end{equation}
Therefore taking into account that $\Phi_{n,s,t}(x)=v_{s}e^{i\left(  2\pi
n+t\right)  x}$ , where $v_{s}$ for

$s=1,2,...m$ are the normalized eigenvector of $C,$ we get
\begin{equation}
\sum\limits_{i=1,2,...,m}(\sum_{n:n\notin A(k,t)}\left\vert \left(
\Psi_{k,j,t},\Phi_{n,i,t}\right)  \right\vert ^{2})=O(\frac{(\ln|k|)^{2}%
}{k^{2}}).
\end{equation}
Since $\left\{  \Phi_{n,i,t}(x):i=1,2,...,m\text{; }n\in Z\right\}  $ is an
orthonormal basis in $L_{2}^{m}\left[  0,1\right]  ,$ we have%
\begin{equation}
\sum\limits_{i=1,2,...,m;n\in A(k,t)}\left\vert \left(  \Psi_{k,j,t}%
(x),\Phi_{n,i,t}(x)\right)  \right\vert ^{2}=1+O(\frac{(\ln|k|)^{2}}{k^{2}})
\end{equation}
Now taking into account that the number of the orthonormal eigenfunctions
$\Phi_{n,i,t}(x)$ for $n\in A(k,t),i=1,2,...,m$ less than $2m$ we obtain the
proof of (29). Instead of (9) using (10), taking into account that the
eigenfunctions of the operator $L_{t}(Q)$ form an orthonormal basis in
$L_{2}^{m}(0,1),$ and arguing as in the proof of (33) we get
\[
\sum\limits_{i=1,2,...,m}(\sum_{n:n\notin A(k,t)}\left\vert \left(
\varphi_{k,j,t},\Psi_{n,i,t}\right)  \right\vert ^{2})=O(\frac{(\ln|k|)^{2}%
}{k^{2}})
\]%
\begin{equation}
\sum\limits_{i=1,2,...,m}(\sum_{n:n\notin A(k,t)}\left\vert \left(
\Phi_{k,j,t},\Psi_{n,i,t}\right)  \right\vert ^{2})=O(\frac{(\ln|k|)^{2}%
}{k^{2}}),
\end{equation}%
\begin{equation}
\sum\limits_{i=1,2,...,m}(\sum_{n:n\in A(k,t)}\left\vert \left(  \Phi
_{k,j,t},\Psi_{n,i,t}\right)  \right\vert ^{2})=1+O(\frac{(\ln|k|)^{2}}{k^{2}%
})
\end{equation}
and the proof of (30)
\end{proof}

\begin{theorem}
\bigskip$(a)$ All big eigenvalues of $L_{t}\left(  Q\left(  x\right)  \right)
$ lie in $\varepsilon_{k}\equiv c_{8}(\mid\frac{\ln|k|}{k}\mid+b_{k})$
neighborhood $U(\varepsilon_{k},\mu_{k,i})$ of the eigenvalues $\mu
_{k,i}(t)=\left(  2\pi k+t\right)  ^{2}+\mu_{i}$ for $\mid$ $k\mid\geq N$,
$i=1,2,...,m$ of $L_{t}(C).$ Moreover, for each big eigenvalues $\mu_{k,j}(t)$
of $L_{t}\left(  C\right)  $ there exists an eigenvalue of $L_{t}(Q(x))$ lying
in $U(\varepsilon_{k},\mu_{k,j}).$

$(b)$ If $t\neq0,\pi$ then the eigenvalues of $L_{t}(Q)$ consist of $m$
sequences (3) satisfying
\begin{equation}
\lambda_{k,j}(t)=\left(  2\pi k+t\right)  ^{2}+\mu_{j}+O(\frac{(\ln|k|)}{k}),
\end{equation}
where $\mu_{1}\leq\mu_{2}\leq...\leq\mu_{m}$ are the eigenvalues of $C,$ and
any normalized eigenfunction $\Psi_{k,j,t}(x)$ corresponding to $\lambda
_{k,j}$ obey
\begin{equation}
\parallel\Psi_{k,j,t}(x)-P\Psi_{k,j,t}(x)\parallel=O(\frac{(\ln|k|)}{k}),
\end{equation}
where $P$ is the orthogonal projection onto the eigenspace corresponding to
$\mu_{k,j}$.

$(c)$ If $\mu_{j}$ is a simple eigenvalue of the matrix $C$ and $t\notin
B(\alpha_{k},\mu_{j}),$ where
\begin{equation}
B(\alpha_{k},\mu_{j})=%
{\textstyle\bigcup\limits_{\substack{n=0,1;\\p=0,1,...,m}}}
\{\pi n+\frac{\mu_{p}-\mu_{j}-\alpha_{k}}{4\pi(2k+n)},\pi n+\frac{\mu_{p}%
-\mu_{j}+\alpha_{k}}{4\pi(2k+n)}\},
\end{equation}
$k\geq N$ and $\alpha_{k}=\sqrt{\varepsilon_{k}},$ then there exists a unique
eigenvalue, denoted by $\lambda_{k,j}(t),$ of $L_{t}\left(  Q\left(  x\right)
\right)  $ lying in $U(\varepsilon_{k},\mu_{k,j}).$ This eigenvalue is a
simple eigenvalue of $L_{t}\left(  Q\left(  x\right)  \right)  $ and the
corresponding normalized eigenfunction $\Psi_{k,j,t}(x)$ satisfy
\begin{equation}
\Psi_{k,j,t}(x)=v_{j}e^{i\left(  2\pi k+t\right)  x}+O(\alpha_{k}),
\end{equation}
where $v_{j}$ is the normalized eigenvector of $C$ corresponding to the
eigenvalue $\mu_{j}.$

$(d)$ Suppose that there exists at least three simple eigenvalues $\mu_{j_{1}%
}<\mu_{j_{2}}<\mu_{j_{3}}$ of $C$ and the set $\{\mu_{j_{p}}+\mu_{i_{p}%
}:p=1,2,3\}$ contains two points for each triple $(i_{1},i_{2},i_{3}),$ where
$i_{1},i_{2},i_{3}=1,2,...,m.$ Then there exist a number $H$ such that
$(H,\infty)\in S(L),$ that is, the number of the gaps in the spectrum of
$L(Q)$ is finite.
\end{theorem}

\begin{proof}
$(a)$ Dividing both sides of (12) by $\left(  \Psi_{k,j,t}(x),\Phi
_{n,i,t}(x)\right)  $ and using (23), (29) we get the proof of the first
statement of $(a).$ Using (30), in the same way, we obtain the proof of the
second statement of $(a).$

$(b).$ If $t\neq0,\pi$ then there exist $N$ such that $t\in T(k)$ for $k>N.$
Therefore instead of (23) using (25), taking into account that $A(k,t)=\{k\}$
for $t\in T(k),$ and repeating the proof of $(a)$ we see that the eigenvalues
$\lambda_{k,j}$ for $j=1,2,...,m$ of $L_{t}\left(  Q\left(  x\right)  \right)
$ lie in $U(\delta_{k},\mu_{k,i})$ for $i=1,2,...,m$ and for each eigenvalue
$\mu_{k,i}$ of $L_{t}(C)$ there exists an eigenvalue of $L_{t}(Q(x))$ lying in
$U(\delta_{k},\mu_{k,i}),$ where $\mid k\mid\geq N,$ $\delta_{k}=c_{9}%
\frac{\ln\mid k\mid}{\mid k\mid}.$

Now we prove that if the multiplicity of the eigenvalue $\mu_{j}$ of $C$ is
$q$ then there exist precisely $q$ eigenvalues of the operator $L_{t}\left(
Q\right)  $ \ lying in $U(\delta_{k},\mu_{k,j})$ for $\mid k\mid\geq N$ . The
eigenvalues of $L_{t}(Q)$ and $L_{t}(C)$ can be numbered in the following way
$\lambda_{k,1}\leq\lambda_{k,2}\leq...\leq\lambda_{k,m}$ and $\mu_{k,1}\leq
\mu_{k,2}\leq...\leq\mu_{k,m}.$ If the matrix $C$ has $r$ different
eigenvalues $\mu_{j_{1}},\mu_{j_{2}},...,\mu_{j_{r}}$ , where $j_{1}%
<j_{2}<...<j_{r}=m,$ with multiplicities $j_{1},$ $j_{2}-j_{1},...,j_{r}%
-j_{r-1}$ then we have%
\begin{align}
\mu_{j_{1}}  &  <\mu_{j_{2}}<...<\mu_{j_{r}};\text{ }\mu_{1}=\mu_{2}%
=...=\mu_{j_{1}};\\
\mu_{j_{1}+1}  &  =\mu_{j_{1}+2}=...=\mu_{j_{2}};...;\mu_{j_{r-1}+1}%
=\mu_{j_{r-2}+2}=...=\mu_{j_{r}}.\nonumber
\end{align}
Suppose that there exist precisely $s_{1},s_{2},...,s_{r}$ eigenvalues of
$L_{t}\left(  Q\left(  x\right)  \right)  $ lying in neighborhoods
$U(\delta_{k},\mu_{k,j_{1}}),U(\delta_{k},\mu_{k,j_{2}}),...,U(\delta_{k}%
,\mu_{k,j_{r}})$ respectively. Since

$\delta_{k}<$ $\frac{1}{2}\min_{p=1,2,...,r-1}\mid\mu_{j_{p}+1}-\mu_{j_{p}%
}\mid$ for $\mid k\mid\geq N$ these neighborhoods are pairwise disjoints. Thus
in the first paragraph of the proof of $(b)$ we have proved that%
\begin{equation}
s_{1}+s_{2}+...+s_{r}=m
\end{equation}
Now we prove that $s_{1}=j_{1},$ $s_{2}=$ $j_{2}-j_{1},$ ..., $s_{r}%
=j_{r}-j_{r-1}.$ Due to the notations the eigenvalues $\lambda_{k,1}%
,\lambda_{k,2}...,\lambda_{k,s_{1}}$ of the operator $L_{t}\left(  Q\right)  $
\ lie in $U(\delta_{k},\mu_{k,1})$ and%
\[
\mid\lambda_{k,j}-\mu_{k,i}\mid>\frac{1}{2}\min_{p:p>j_{1}}\mid\mu_{1}-\mu
_{p}\mid
\]
for $j\leq s_{1}$ and $i>j_{1}.$ Therefore it follows from
\begin{equation}
(\lambda_{k,j}(t)-\mu_{k,i})(\Psi_{k,j,t},\Phi_{k,i,t})=((Q(x)-C)\Psi
_{k,j,t}(x),\Phi_{k,i,t}(x))
\end{equation}
and (25) that%
\begin{equation}
\sum_{i>j_{1}}\left\vert \left(  \Psi_{k,j,t},\Phi_{k,i,t}\right)  \right\vert
^{2})=O(\frac{(\ln|k|)^{2}}{k^{2}}),\forall j\leq s_{1}.
\end{equation}
Using this and (33), and taking into account that $A(k,t)=\{k\}$ for $t\in
T(k)$ ( see (6)), we conclude that there exists normalized eigenfunction
corresponding to $\mu_{k,1}=\mu_{k,2}=...=\mu_{k,j_{1}}$ and denoted by
$\Phi_{k,j,t}(x)$ such that%
\begin{equation}
\Psi_{k,j,t}(x)=\Phi_{k,j,t}(x)+O(\frac{\ln|k|}{k})
\end{equation}
for $j\leq s_{1}.$ Since $\Psi_{k,1,t}(x),\Psi_{k,2,t}(x),...,\Psi_{k,s_{1}%
,t}(x)$ are orthonormal system we have
\[
(\Phi_{k,j,t}(x),\Phi_{k,i,t}(x))=\delta_{i,j}+O(\frac{\ln|k|}{k})
\]
for $i,j=1,2,...,s_{1}.$ This formula imply that the dimension $j_{1}$ of the
eigenspace of $L_{t}(C)$ corresponding to the eigenvalue $\mu_{k,1}=$ $\mu
_{1}+(2\pi k+t)^{2}$ is not less than $s_{1}.$ Thus $s_{1}\leq j_{1}.$ In the
same way we prove that $s_{2}\leq$ $j_{2}-j_{1},$ ..., $s_{r}\leq
j_{r}-j_{r-1}.$ These inequalities with the equalities (41), $j_{r}=m$ imply
that $s_{1}=j_{1},$

$s_{2}=$ $j_{2}-j_{1},$ ..., $s_{r}=j_{r}-j_{r-1}.$ Therefore taking into
account that ( see (3),(4)) the eigenvalues of the operator $L_{t}\left(
Q\right)  $ consist of $m$ sequences satisfying (4) we get the prove of (36).
The proof of (37) follows from (44).

$(c)$ To consider the simplicity of $\mu_{k,j}(t)$ and $\lambda_{k,j}(t)$ for
$k\geq N$ we introduce the set
\begin{equation}
S(k,j,n,i)=\{t\in\lbrack-\frac{\pi}{2},\frac{3\pi}{2}):\mid\mu_{k,j}%
(t)-\mu_{n,i}(t)\mid<\alpha_{k}\}
\end{equation}
for $(n,i)\neq(k,j).$ It follows from (5) that $S(k,j,n,i)=\emptyset$ for
$n\neq k,-k,-k-1.$ If $\mu_{j}$ is a simple eigenvalue of the matrix $C$ then
$S(k,j,k,i)=\emptyset$ for $i\neq j.$ It remains to consider the sets
$S(k,j,-k,i),$ $S(k,j,-k-1,i)$ By direct calculations of the differences
$\mu_{k,j}(t)-\mu_{-k,i}(t),$ $\mu_{k,j}(t)-\mu_{-k-1,i}(t)$ one can easily
verify that
\begin{equation}
S(k,j,-k,i)=(\frac{\mu_{i}-\mu_{j}-\alpha_{k}}{8\pi k},\frac{\mu_{i}-\mu
_{j}+\alpha_{k}}{8\pi k})
\end{equation}%
\begin{equation}
S(k,j,-k-1,i)=(\pi+\frac{\mu_{i}-\mu_{j}-\alpha_{k}}{4\pi(2k+1)},\pi+\frac
{\mu_{i}-\mu_{j}+\alpha_{k}}{4\pi(2k+1)})
\end{equation}
These relations with (38) imply that
\begin{equation}
B(\alpha_{k},\mu_{j})=%
{\textstyle\bigcup\limits_{\substack{n\in\mathbb{Z};(n,i)\neq
(k,j)\\i=0,1,...,m}}}
S(k,j,n,i)=%
{\textstyle\bigcup\limits_{\substack{n=-k,-k-1\\i=0,1,...,m}}}
S(k,j,n,i)
\end{equation}
Therefore it follows from the definition of $S(k,j,n,i)$ ( see (45)) that if

$\ t\notin B(\alpha_{k},\mu_{j}),$ then
\begin{equation}
\mid\mu_{k,j}(t)-\mu_{n,i}(t)\mid\geq\alpha_{k}%
\end{equation}
for all $(n,i)\neq$ $(k,j).$ Hence $\mu_{k,j}(t)$ is a simple eigenvalue for
$t\notin B(\alpha_{k},\mu_{j}).$ By $(a)$ there exist an eigenvalue of
$L_{t}(Q)$ lying in $U(\varepsilon_{k},\mu_{k,j}).$ Denote this eigenvalue by
$\lambda_{k,j}(t)$. Thus $\lambda_{k,j}$ is an eigenvalue of $L_{t}(Q(x))$
satisfying
\begin{equation}
\mid\lambda_{k,j}(t)-\mu_{k,j}(t)\mid<\varepsilon_{k}.
\end{equation}
It follows from this and (49) that
\begin{equation}
\mid\lambda_{k,j}(t)-\mu_{n,i}(t)\mid>\frac{1}{2}\alpha_{k}.
\end{equation}
for $n\in A(k,t),$ $i=1,2,...,m,$ and $(n,i)\neq(k,j)$ since $k\geq N,$
$\alpha_{k}=\sqrt{\varepsilon_{k}}$ and $\varepsilon_{p}\rightarrow0$ as
$p\rightarrow\infty.$ Let $\Psi_{k,j,t}(x)$ be any normalized eigenfunction
corresponding to $\lambda_{k,j}.$ Using (12), (23), definition of
$\varepsilon_{k},$ and (51), we obtain%
\[
(\Psi_{k,j,t}(x),\Phi_{n,i,t}(x))=O\left(  \alpha_{k}\right)
\]
for $(n,i)\neq(k,j)$ and $n\in A(k,t).$ This and (35) imply that $\Psi
_{k,j,t}(x)$ satisfies (39). Thus we have proved that (39) holds for any
normalized eigenfunction corresponding to any eigenvalue lying in
$U(\varepsilon_{k},\mu_{k,j})$. If there are two different eigenvalue or a
multiple eigenvalue lying in $U(\varepsilon_{k},\mu_{k,j})$ then there are two
orthonormal eigenfunction satisfying (39) which is impossible. Hence there
exists unique eigenvalue $\lambda_{k,j}(x)$ lying in $U(\varepsilon_{k}%
,\mu_{k,j})$ and this is a simple eigenvalue.

$(d)$ In $(c)$ we proved that if $k\geq N$ , $t\notin B(\alpha_{k},\mu_{j_{p}%
})$, where $p=1,2,3,$ then there exists a unique eigenvalue of $L_{t}(Q)$
lying in $U(\varepsilon_{k},\mu_{k,j_{p}}).$ We denoted this unique eigenvalue
by $\lambda_{k,j_{p}}(t)$ and proved that it it a simple eigenvalue. Let us
prove that $\lambda_{k,j_{p}}(t)$ is a continuous function at $t_{0}\in
\lbrack-\frac{\pi}{2},\frac{3\pi}{2})\backslash B(\alpha_{k},\mu_{j_{p}}).$
Since $\lambda_{k,j_{p}}(t_{0})$ is a simple eigenvalue it is a simple root of
the characteristic determinant $\Delta(\lambda,t)$ of the operator $L_{t}(Q).$
Therefore there exists a neighborhood $U(t_{0})$ of $t_{0}$ and a continuous
in $U(t_{0})$ function $\Lambda(t)$ such that $\Lambda(t_{0})=\lambda
_{k,j_{p}}(t_{0}),$ $\Lambda(t)$ is an eigenvalue of $L_{t}(Q)$ for $t\in
U(t_{0})$ and
\begin{equation}
\mid\Lambda(t)-\mu_{k,j_{p}}(t)\mid<\varepsilon_{k},\forall t\in U(t_{0}),
\end{equation}
since $\mid\Lambda(t_{0})-\mu_{k,j_{p}}(t_{0})\mid<\varepsilon_{k}$ and the
functions $\Lambda(t),$ $\mu_{k,j_{p}}(t)$ are continuous. Now taking into
account that there exists unique eigenvalue of $L_{t}(Q)$ lying in
$U(\varepsilon_{k},\mu_{k,j_{p}}),$ we obtain that $\Lambda(t)=\lambda
_{k,j_{p}}(t)$ for $t\in U(t_{0}),$ and hence $\lambda_{k,j_{p}}(t)$ is
continuous at $t_{0}.$ Now we prove that there exists $H$ such that
\begin{equation}
(H,\infty)\subset\{\lambda_{k,j_{p}}(t):t\in\lbrack-\frac{\pi}{2},\frac{3\pi
}{2})\backslash B(\alpha_{k},\mu_{j_{p}}),\text{ }k=N,N+1,...,\}.
\end{equation}
It is clear that
\begin{equation}
(h,\infty)\subset\{\mu_{k,j_{p}}(t):t\in\lbrack-\frac{\pi}{2},\frac{3\pi}%
{2}),k=N,N+1,...,\},
\end{equation}
where $h=\mu_{N,,j_{3}}(-\frac{\pi}{2}).$ Since $\mu_{k,j_{p}}(t)=(2\pi
k+t)^{2}+\mu_{j_{p}}$ is increasing function for $k\geq N,$ it follows from
(48), (47), (46) that
\[
\{\mu_{k,j_{p}}(t):t\in B(\alpha_{k},\mu_{j_{p}})\}\subset%
{\textstyle\bigcup\limits_{\substack{n=0,1;\\i=0,1,...,m}}}
C(k,j_{p},i,\alpha_{k},n),
\]
where $C(k,j_{p},i,\alpha_{k},n)=\{x\in\mathbb{R}:\mid x-(\pi(2k+n))^{2}%
-\frac{\mu_{i}+\mu_{j_{p}}}{2}\mid<\alpha_{k}\}.$ Therefore
\[
\{\mu_{k,j_{p}}(t):t\in\lbrack-\frac{\pi}{2},\frac{3\pi}{2})\backslash
B(\alpha_{k},\mu_{j_{p}}),k\geq N\}\supset(h,\infty)\backslash%
{\textstyle\bigcup\limits_{\substack{n=0,1;k\geq N\\i=0,1,...,m}}}
C(k,j_{p},i,\alpha_{k},n)
\]
Now using (50) and the continuity of $\lambda_{k,j_{p}}(t)$ on $[-\frac{\pi
}{2},\frac{3\pi}{2})\backslash B(\alpha_{k},\mu_{j_{p}})$ we get
\[
\{\lambda_{k,j_{p}}(t):t\in\lbrack-\frac{\pi}{2},\frac{3\pi}{2})\backslash
B(\alpha_{k},\mu_{j_{p}}),k\geq N\}\supset(H,\infty)\backslash%
{\textstyle\bigcup\limits_{\substack{n=0,1;k\geq N\\i=0,1,...,m}}}
C(k,j_{p},i,2\alpha_{k},n),
\]
{}

where $H=h+1.$ Thus
\[%
{\textstyle\bigcup\limits_{p=1,2,3}}
((H,\infty)\backslash%
{\textstyle\bigcup\limits_{n=0,1;k\geq N;i=0,1,...,m}}
C(k,j_{p},i,2\alpha_{k},n))\subset S(L)
\]
To prove that $(H,\infty)\subset S(L)$ it is enough to show that the set
\[%
{\textstyle\bigcap\limits_{p=1,2,3}}
(%
{\textstyle\bigcup\limits_{n=0,1;k\geq N;i=0,1,...,m}}
C(k,j_{p},i,2\alpha_{k},n))
\]
is empty. If this set contains an element $x,$ then%
\[
x\in%
{\textstyle\bigcup\limits_{n=0,1;k\geq N;i=0,1,...,m}}
C(k,j_{p},i,2\alpha_{k},n)
\]
for all $p=1,2,3.$ Using this and the definition of $C(k,j_{p},i,2\alpha
_{k},n),$ \ we obtain that there exist $k\geq N;$ $n=0,1$ and $i_{p}$ such
that
\[
\mid x-(\pi(2k+n))^{2}-\frac{\mu_{j_{p}}+\mu_{i_{p}}}{2}\mid<2\alpha_{k}%
\]
for all $p=1,2,3$ and hence $\mid\frac{\mu_{j_{q}}+\mu_{i_{q}}}{2}-\frac
{\mu_{j_{p}}+\mu_{i_{p}}}{2}\mid<4\alpha_{k}$ for all $p,q=1,2,3.$ Since
$\alpha_{k}\rightarrow0$ as $k\rightarrow\infty$ the last inequality imply
that $\mu_{j_{1}}+\mu_{i_{1}}=\mu_{j_{2}}+\mu_{i_{2}}=\mu_{j_{3}}+\mu_{i_{3}}$
which contradicts the condition of $(d).$
\end{proof}

\end{document}